\documentclass{svmult}       
\usepackage[latin1]{inputenc}
\usepackage{makeidx}         
\usepackage{graphicx}        
\usepackage{multicol}        
\usepackage[bottom]{footmisc}

\usepackage{booktabs}      
\usepackage{dcolumn}       
\usepackage{url}           
\usepackage{amssymb}
\usepackage{amsmath}
\usepackage{listings}      
\lstdefinelanguage[rngts_keywords]{xml}[]{xml}{%
  morekeywords={RNG_TEST_SUITE_RESULT, RNG, SEED, TEST, PARAMETERS, PARAMETER,
  ANALYZE, CHI_SQUARE, KOLMOGOROV_SMIRNOV, RESULTS, PASSED, FAILED, RESULT}%
}                          
\usepackage{mdwlist}       
\usepackage[numbers]{natbib}    

\makeatletter \def\@captype{table}\makeatother

\makeatletter
\DeclareRobustCommand{\Cpp}
{\valign{\vfil\hbox{##}\vfil\cr
   \textsf{C\kern-.1em}\cr
   $\hbox{\fontsize{\sf@size}{0}\textbf{+\kern-0.05em+}}$\cr}%
}
\makeatother

\makeindex             

\begin{document}

\title*{A Generic Random Number Generator Test Suite}

\author{Mario Rütti \inst{1} \and Matthias Troyer \inst{2} \and Wesley
  P. Petersen \inst{1}}

\institute{Seminar for Applied Mathematics, ETH Z\"urich, CH-8092 Z\"urich
\url{maruetti@comp-phys.org}
\and Institute for Theoretical Physics, ETH Z\"urich, CH--8093 Z\"urich }

\maketitle

\begin{abstract}
The heart of every \index{Monte Carlo}Monte Carlo simulation is a
source of high quality random numbers and the generator has to be
picked carefully. Since the ``Ferrenberg affair''
\cite{Ferrenberg:1992} it is known to a broad community that statistical tests 
alone do not suffice to determine the quality  of a generator, but also application-based 
tests are needed. With the inclusion of an extensible random number library and the 
definition of a generic interface into the revised \Cpp{} standard \cite{cpp} it will 
be important to have access to an extensive \Cpp{} random number test suite. Most currently 
available test suites are limited to a subset of tests are written in Fortran or C and 
cannot easily be used with the \Cpp{} random number generator library.

In this paper we will present a generic random number test suite written in \Cpp. The 
suite is based on the Boost reference implementation \cite{boost} of the forthcoming 
\Cpp{} standard random number generator library. The Boost implementation so far contains 
most modern random number generators. Employing generic programming techniques the test 
suite is flexible, easily extensible and can be used with any random number generator 
library, including those written in C and Fortran. Test results are produced in an XML 
format, which through the use of XSLT transformations allows extraction of summaries or 
detailed reports, and conversion to HTML, PDF, PostScript or any other format. 

At this time, the test suite contains a wide range of different tests, including the 
standard tests described by Knuth \cite{Knuth:1997}, Vattulainen's physical tests 
\cite{Vattulainen:1994}, parts of Marsaglia's Diehard \cite{Marsaglia:1985} test suite,
and a number of number of newer test such as Gonnet's repetition test
\cite{Gonnet-repeating-test-03}
\end{abstract}

\section{Introduction}
\label{sec:introduction-1}
Testing random number generators is a really important part in the selection process for
an adequate \index{random number generator}random number generator. But testing random number generators is also a very
serious problem and not so easy to manage. In this paper we will present a generic \index{framework}framework
(random number generator test suite,
\index{\textsc{RNGTS}}\textsc{RNGTS}) to manage, perform and analyze 
random number generator tests. The framework
is written in \Cpp{} and makes an extensive use of the \index{generic programming technique}generic programming technique. This 
technique allows the usage of different data types and supports all different random number 
generators which fulfill the specification in the upcoming \Cpp{} standard for random number
generators. An other important point is the flexibility of the framework. Our attempt offers
a flexible and simple extensible solution for further tests and other generators.
Using \index{XML}XML to format the output opens a wide field of standardized transformations to other
formats. 

\subsection{Random Number Generator Tests}
\label{sec:rand-numb-gener-1}
Using a \index{random number generator}random number generator is one
thing, the knowledge about its \index{performance}performance another.
Testing the performance means first thinking about the expected properties of the random
numbers, for example its application specific requirements. In a second step the most 
suitable tests have to be chosen. This selection is not easy and often people make 
``overall'' tests of the generator. The problem is that we probably know what the test does,
but one rarely knows which one the sensitive parameters in the application are.
So, it is preferable to run more tests than fewer.\par
The best known source is Knuth's book \cite{Knuth:1997}, but there are many other tests
available, some included in test benches, other as single procedures. The best known test bench
is probably the \textsc{Diehard} test suite from Marsaglia \cite{diehard:old} of which
a new version has been released \cite{diehard:new}. An other \index{test bench}test bench, more specialized on 
\index{physical tests}physical tests, is available from Vattulainen \cite{Vattulainen:1999}. Furthermore there is 
the \textsc{SPRNG} library \cite{sprng} for parallel generators. In addition to the numerated 
test benches there are a lot of stand alone procedures like ``Maurers Universal Test'' 
\cite{Maurer90b} or the ``Repeating Time Test'' \cite{Gonnet-repeating-test-03}. Different 
tests are a dime a dozen, so we collected the most popular ones in table \ref{tab:avail-tests}. 
All tests in the table having a ``Class Name'' are available in the \textsc{RNGTS} framework.

\subsection{Generic Programming}
\label{sec:generic-programming}
\index{Generic programming}Generic programming in \Cpp{} became popular with the introduction of the \index{Standard Template Library
(STL)}Standard Template Library
(STL) into the \Cpp{} Standard. The STL is a collection of often used Containers and Algorithms.
The principle of generic programming is simple, instead of fixing the view onto particular
data structures, the main interest are lying on the algorithms. This requires the introduction
of a special data type handling mechanism. In \Cpp{} this technique is called template programming.\par
One simple example is the maximum operator, the ordinary maximum function. In an 
ordinary way, the maximum operator may be defined as follows.
\begin{lstlisting}[texcl]{}
const int max(const int& a, const int& b)
{ return (a > b) ? a : b; }
\end{lstlisting}
This method works well with integers, but we get into trouble if there are two doubles to
compare. One solution is the replacement of all int's with double's. But the next problem
occurs if the values are of type byte. The solution for this rapidly multiplying dilemma is, 
as mentioned above, templates. The specific data type is replaced by a generic placeholder 
which must be defined at compile time. Now, the maximum operator can be defined as followed.
\begin{lstlisting}[texcl]{}
template <typename T>
const T max(const T& a, const T& b)
{ return (a > b) ? a : b; }
\end{lstlisting}
The precondition for making this method work is that the comparison operator ``\textgreater'' must be
defined.\par
The diversity of the STL gives an idea how powerful the usage of the template concept is.
But this is only the tip of the iceberg, one can do much more. There is a technique called
``\index{Template Metaprogramming}Template Metaprogramming'', based on
templates. This syntax allows \index{compile time programming}compile
time programming,
e.~g. nifty things like compile time if's or decisions about used data types. A number of different
applications of this template technique can be found in the \textsc{Boost} library \cite{boost}
as ``Generic Programming'' and ``Template Metaprogramming''.

\subsection{The Boost library}
\label{sec:boost-library}
This section gives a short description of the \index{\textsc{Boost}}\textsc{Boost} project and presents the 
``\textsc{Boost} Random Number Library''.\par
The following quotation from the \textsc{Boost} home page describes the project in 
a few words.
\begin{quotation}
  The Boost web site provides free peer-reviewed portable \Cpp{} source libraries. The emphasis 
is on libraries which work well with the \Cpp{} Standard Library. The libraries are intended to 
be widely useful, and are in regular use by thousands of programmers across a broad spectrum 
of applications.\\
A further goal is to establish "existing practice" and provide reference implementations so 
that Boost libraries are suitable for eventual standardization. Ten Boost libraries will be 
included in the \Cpp{} Standards Committee's upcoming \Cpp{} Standard Library Technical Report as 
a step toward becoming part of a future \Cpp{} Standard.\\
\end{quotation}

The project contains already a huge collection of libraries and one of them is the \texttt{Random} 
library. The \texttt{Random} library is one of the libraries included in the upcoming \Cpp{} Standard. 
Working with this random library means working with the future of \Cpp.

The \texttt{Random} library assembles a specific random number generator from to parts.
One part is an ``engine'' (the raw random number generator), the other part is a ``distribution''. 
This two parts together are represented by a so called \texttt{variate\_generator}. For both 
parts there are a number of predefined types. Some examples are shown in table \ref{tab:generators} 
and \ref{tab:distributions}.

\begin{table}[htbp]
  \centering
  \begin{tabular}{ll}
    \toprule
    Boost Random Number Engine    & Example               \\
    \cmidrule(lr){1-2}
    \texttt{linear\_congruential} & GGL                   \\
    \texttt{additive\_combine}    & L'Ecuyer 1988         \\
    \texttt{lagged\_fibonacci}    & R1279                 \\
    \texttt{mersenne\_twister}    & MT19937               \\
    \texttt{shuffle\_output}      & Bays-Durham shuffle   \\
    \ldots                                                \\
    \bottomrule
  \end{tabular}
  \caption{Some \textsc{Boost} random library generators}
  \label{tab:generators}
\end{table}

\begin{table}[htbp]
  \centering
  \begin{tabular}{ll}
    \toprule
    Boost Random Number Distribution     & Example              \\
    \cmidrule(lr){1-2}
    \texttt{uniform\_01}                 & $x \in \mathbb{R}, x \in [0,1)$ \\
    \texttt{uniform\_int}                & $x \in \mathbb{N}$ \\
    \texttt{uniform\_real}               & $x \in \mathbb{R}, \quad x \in [\mathrm{min},\mathrm{max})$ \\
    \texttt{uniform\_on\_sphere}         & $x \in \mathbb{S}^n$ \\
    \texttt{normal\_distribution}        & $P(x) = \frac{1}{\sigma \sqrt{2\pi}} e^{-(x-\mu)^2/(2\sigma^2)}$ \\
    \texttt{exponential\_distribution}   & $P(x) = \lambda e^{-\lambda x}$ \\
    \ldots & \\
    \bottomrule
  \end{tabular}
  \caption{Some \textsc{Boost} random library distributions}
  \label{tab:distributions}
\end{table}

To manufacture a desired generator, one has only to specify the desired components and
put them together. The code snippet below shows an example.
\begin{lstlisting}[texcl]{}
#include <boost/random.hpp>         // the required libraries

boost::mt19937 rng;                 // the engine is a Mersenne Twister
boost::uniform_int<> six(1,6)       // the distribution maps to integer 1..6

// glue the components together
boost::variate_generator<boost::mt19937, boost::uniform_int<> >  
    die(rng, six);             

int x = die();                      // use the generator...
\end{lstlisting}

\section{Random Number Generator Test Suite Framework}
\label{sec:rand-numb-gener-2}

\subsection{Statistics}
\label{sec:statistics}
After running a test one has to check if the obtained result is the expected one or something
different. To classify the differences there are a variety of different \index{statistical methods}statistical methods.
There are three different statistical methods most frequently used. These are \emph{Chi-Square} 
($\chi^2$), \emph{Kolmogorov-Smirnov} and \emph{Gaussian} test. To give an introduction, 
there is a short description below. As the final output of each statistic, one or more 
probabilities in the $[0,1]$ range are expected. These values can provide a \index{confidence level}confidence level.

\begin{description}
\item[\index{Chi-Square ($\chi^2$)}Chi-Square ($\chi^2$)] 
  The $\chi^{2}$-Test is perhaps the best known statistical test. It is based on a comparison 
  between the empirical data and a theoretical distribution. The empirical data are the results
  of the random process.
  \begin{equation}
    \label{eq:1}
    \chi^{2} = \sum_{i=1}^{k}\frac{(n_{i}-np_{i})^{2}}{np_{i}} 
    = \frac{1}{n} \sum_{i=1}^{k}\left(\frac{n_{i}^2}{p_{i}}\right)-n
  \end{equation}

\item[\index{Kolmogorov-Smirnov}Kolmogorov-Smirnov] 
  The $\chi^{2}$ test is applied when observations can fall into a finite number of categories. 
  But normally one will consider random quantities which may assume an infinite number of values.
  In this test, the random number generators distribution function $F_{n}(x)$ is compared to the
  theoretical cumulative distribution function$F(x)$. 
  \begin{eqnarray}
    K_{n}^{+} &=& \sqrt{n} \max_{-\infty < x < \infty} \left(F_{n}(x) -
      F(x) \right) =
    \sqrt{n} \max_{1 < i < n} \left(\frac{i}{n} - F(X_{i})\right) \\
    K_{n}^{-} &=& \sqrt{n} \max_{-\infty < x < \infty} \left(F(x) -
      F_{n}(x) \right) =
    \sqrt{n} \max_{1 < i < n} \left(F(X_{i}) - \frac{i-1}{n}\right)
  \end{eqnarray}

\item[\index{Gaussian}Gaussian] 
  The Gaussian test is a little different from the $\chi^2$ or the Kolmogorov-Smirnov test. 
In these two tests the expected distribution function is compared with the measured distribution 
function and based on the difference some indicators are calculated. In the Gaussian test a physical 
view is used. If a measurement is done, it is known that, even if the best tools are used, the result
depends on a number of ruleless and uncontrolled parameters. These measurement errors are random and
a combination of different single errors.
This statistic calculates the probability for the deviation from the mean value. The used formula
is 
\begin{equation}
  \mathrm{p(x)} = \frac{1}{2} - \frac{1}{2} \mathrm{erf} \left(\frac{1}{\sqrt{2}}x\right)
\end{equation}
where $\mathrm{erf}$ denotes the ``error function'' 
$\mathrm{erf(z)} = \frac{2}{\sqrt{\pi}} \int_{0}^{z} e^{-t^2} dt$.
\end{description}

The \textsc{RNGTS}-Framework implements these three statistics as base classes for the tests. They
are named as \texttt{chisquare\_test}, \texttt{ks\_test} and \texttt{gaussian\_test} and implement
the needed statistical methods. Using one of these classes is not mandatory, but recommended. A
benefit is the possibility of using other helpers like \texttt{iterate\_test}.

\subsection{Tests}
\label{sec:tests}
Well known sources for random number generator tests are the book of Knuth \cite{Knuth:1997} and 
Marsaglias \textsc{Diehard} test suite \cite{diehard:old} and \cite{diehard:new}. But there are 
many more tests, perhaps not so nicely packaged as in the works mentioned above, but still useful. 
Each of these test checks a specific property of the tested generator, e.~g. the uniformity of 
binned values or the uniformity of 'continuous' values. The most popular tests are shown in table 
\ref{tab:avail-tests}. The same table contains also the name of the class which implements the 
test, if it is already implemented.\par
To make a test work with the \textsc{RNGTS} framework, one only has to implement a specific
interface, or in other words, there are some methods in a test class which are needed by the
framework. The implementation may be done from scratch or if one of the often used statistics
is needed, by deriving from a given base class. 

\paragraph{Implementing a test:}
The following listing shows the base of each test, containing all required methods. This methods
have to be implemented if the test is not derived from one of the provided test base classes.

\begin{lstlisting}[texcl]{}
#include "xml_helper.h"     // XML output functions

class the_new_test 
{
public:                   
  the_new_test(...);

  template< class RNG >
  void run(RNG& rng);

  std::string test_name() const;

  template < class InputIterator >
  void analyze(xml_stream& out, 
               InputIterator cl_begin, 
               InputIterator cl_end) const;

  void print_parameters(xml_stream& out) const;
}
\end{lstlisting}

\begin{description}
\item[constructor] The constructor has to take all parameters needed for further calculations,
  e.g. the number of random numbers to test or the number of test runs.
\item[run] The run method is the core of every test, this method has to run the test sequence
  and calculate the appropriate statistic. The generator to use is passed as \texttt{rng} and
  fulfills the \textsc{Boost} random number generator specification.
\item[test\_name] This method has to return the name of the test which is printed in the output.
\item[analyze] The analyze method checks the test results against a given set of confidence 
  levels, which are passed by iterators. The check also includes the decision if a test passed
  or failed at the confidence level. The output is written onto an XML stream. Several methods 
  for producing valid XML output are provided.
 The output has to fulfill the XML Schema.
\item[print\_parameters] This method has to print all relevant parameters to the given XML
  stream.
\end{description}

\begin{table}[htbp] 
  \centering
  \begin{tabular}{*{2}{l}}
  \toprule
  Test & Class Name \\
  \cmidrule(lr){1-1}\cmidrule(lr){2-2}
  Equidistribution Test (Frequency Test)      & ks\_uniformity\_test       \\
                                              & chisqr\_uniformity\_test   \\ 
  Gap Test                                    & gap\_test                  \\
  Ising Model Test                            & ising\_model\_test         \\
  n-block test                                & n\_block\_test             \\
   
  Serial Test                                 & serial\_test               \\
  Poker Test (Partition Test)                 & poker\_test                \\
  Coupon collector's Test                     & coupon\_collector\_test    \\
  Permutation Test                            & permutation\_test          \\
  Run Test                                    & runs\_test                 \\
  Maximum of t Test                           & max\_of\_t\_test           \\
  Collision Test (Hash Test)                  & collision\_test            \\
  Serial correlation Test                     & serial\_correlation\_test  \\

  Birthday-Spacing's Test                     & birthday\_spacing\_test    \\
  Overlapping Permutations Test               & overlapping\_permutations\_test \\
  Ranks of $31 \times31$ and $32 \times32$ matrices Test &  bin\_rank\_chisqr\_test \\
  Ranks of $6 \times 8$ Matrices Test         & bin\_rank\_ks\_test        \\
  Monkey Tests on 20-bit Words                & monkey\_20bit\_test        \\
  Monkey Tests OPSO,OQSO,DNA                  & monkey\_[OPSO\textbar OQSO\textbar DNA]\_test \\
  Count the 1`s in a Stream of Bytes          & count\_ones\_stream\_test  \\
  Count the 1`s in Specific Bytes             & count\_ones\_bytes\_test   \\
  Parking Lot Test                            & parking\_lot\_test         \\
  Minimum Distance Test                       & minimum\_distance\_test    \\
  Random Spheres Test                         & random\_sphere\_test       \\
  The Sqeeze Test                             & squeeze\_test              \\
  Overlapping Sums Test                       & overlapping\_sums\_test    \\
  The Craps Test                              & craps\_test                \\

  Sum of distributions (for parallel streams) &                            \\
  FFT                                         &                            \\
  Blocking Test                               &                            \\  
  2-d Random Walk                             &  random\_walk\_test        \\
  Random Walkers on a line ($S\_n$ Test)      &                            \\
  2D Intersection Test                        &                            \\
  2D Height Correlation Test                  &  height\_corr2d\_test      \\
  Repeating Time Test                         &  repetition\_test          \\
  Gorilla Test                                &  gorilla\_test             \\
  GCD Test                                    &  gcd\_test                 \\
  Maurers Universal Test                      &  maurers\_universal\_test  \\
  \bottomrule
  \end{tabular}
  \caption{Available tests in the \textsc{RNGTS} framework\label{tab:avail-tests}. 
  A detailed description or/and the original sources may be found in \cite{Ruetti}.}
\end{table}

\subsection{Other Utilities}
\label{sec:other-utilities}
To perform a widespread analysis of a random number generator, the most important using 
is a variate of different tests. Occasionally, people want to test particular behavior 
of the generator, or they have to repeat a test multiple times or for different seeding
strategies, etc. A short overview of the possibilities for included utilities are itemized
below.
\begin{description}
\item[count\_fails\_test] Repeats a test several times and count the number of failings 
  per confidence level
\item[iterate\_test] Repeats a test several times and calculates a Kolmogorov-Smirnov 
  statistic from the results
\item[parallel\_rng\_imitator] Combines multiple random number generators to one generator
\item[rng\_bit\_extract]Interpret a random number as bit field and combine a sequence of 
  successive bits to a new number
\item[rng\_bit\_test] Interpret a random number as bit field and generate sub-fields by 
  shifting a specific mask over the field
\item[rng\_file] Reads random numbers from a binary file
\item[rng\_wrapper] Wrap external generators (e.g. Fortran, C) to use in the framework
\end{description}

\subsection{Results}
\label{sec:results}
To rate a tested random number generator, one needs an overview over all test results
with different seeds at different confidence levels. Additional interesting information
to ensure the reproducibility of performed tests should also be saved.\\
This variety of different information are stored in an \index{XML}XML structure. An XML 
structure allows transformations of the raw results to different target representations
like HTML or \LaTeX. These two transformations are implemented using \index{xslt}xslt (XML Stylesheet 
Language Translation). The ability of newest Internet browsers to format and display XML 
files makes the generation of a report really simple. If a report in \LaTeX{} style is
required, an xslt processor (like xsltproc) has to be used.

\subsection{An Example}
To give an impression of the usage of the \textsc{RNGTS} framework, a short example is 
presented. As we can see here, the usage is intuitive and straight forward. \par

\label{sec:an-example}
\begin{lstlisting}[texcl]{}
#include <boost/random.hpp>               // import the boost random library
#define PRINT_STATUS                      // print status information at run time
#include <rngts/rng_test_suite.h>         // import the rngts framework
#include <rngts/chisqr_uniformity_test.h> // import all required tests

int main()
{
  rng_test_suite<> rngts;                 // create test-suite

  // specify confidence levels
  rngts.add_confidence_level(0.05);       // add  5\% confidence level
  rngts.add_confidence_level(0.95);       // add 95\% confidence level 

  rngts.add_seed(331);                    // add seed 331
  rngts.add_seed(667790);                 // add seed 667790

  // register random number generators to test
  // using Mersenne-Twister and RanLux from boost
  rngts.register_rng<boost::mt19937>("mt-19937");    
  rngts.register_rng<boost::ranlux64_3_01>("rl-64"); 

  // create and initialize a chisquare test object
  chisqr_uniformity_test chi(100000, 256);
  // register the chisquare test
  rngts.register_test<chisqr_uniformity_test>(chi);

  std::ofstream file_out("results.xml");  // file to write output
  try
  {
    rngts.run_test(file_out, true);       // run all tests
  } catch (std::exception& e)            
  {
    std::cout << "failure : " << e.what();
  }
  file_out.close();                       // close output file
}
\end{lstlisting}

If the example above is compiled and executed, the file \texttt{results.xml}
is written. A part of this file is printed below. The extract shows the results
of a $\chi^2$ test of the \emph{Mersenne Twister 19937} generator. The test was
performed with a seed of 331 and no ``warmup'' runs were executed. These are 
the configuration data of the generator section. In the test section, the needed
test parameters are contained in the \texttt{PARAMETERS} block, and the \texttt{ANALYZE} 
block contains the test results. The results are included in the test name tag 
or in generalized result tags. Information whether a test passed or failed is 
represented by tags named \texttt{PASSED} or \texttt{FAILED} which also contain 
the appropriate confidence level.

\begin{lstlisting}[texcl,language={[rngts_keywords]{xml}}]{}
<?xml version="1.0" ?><?xml-stylesheet href="xml2html.xsl" type="text/xsl"?>
<RNG_TEST_SUITE_RESULT  date="2004-04-26">
   <RNG  name="mt-19937"  warmup="0">
      <SEED  seed="331">
         <TEST  name="Chi-Square-Uniformity-Test">
            <PARAMETERS>
               <PARAMETER  name="Number of Numbers"  value="100000"/>
               <PARAMETER  name="Number of Classes"  value="256"/>
            </PARAMETERS>
            <ANALYZE>
               <CHI_SQUARE  chi2="242.33"  probability="0.706"  dof="255">
                  <PASSED  confidenceLevel="0.05"/>
                  <PASSED  confidenceLevel="0.95"/>
               </CHI_SQUARE>
            </ANALYZE>
         </TEST>
      </SEED>
      .....
   </RNG>
</RNG_TEST_SUITE_RESULT>
\end{lstlisting}

\section{Summary}
\label{sec:summary}
The presented \textsc{RNGTS} framework is a powerful but manageable tool to test
the performance of \Cpp{} standard compatible random number generators. It provides
a variety of random number generator tests an some other helpfully utilities. The
ease of transforming XML to other formats allows postprocessing results to different
representation like HTML or \LaTeX. \par
But, one has to keep in mind, the final decision about using a generator or not
can not be done by any program, because the usability strongly depends on the 
application. But, this framework should be a help to simply the decision.

\section{Sources}
The \textsc{RNGTS} framework is available at \url{http://www.comp-phys.org/rngts}.
At this point one can also find a detailed description of the framework, some examples
and related material.

\section{Acknowledgements}
We want to thank to the Swiss Mathematical Society (SMG) for partially support this
work.

\bibliography{mrmtwp}
\bibliographystyle{plainnat}

\printindex
\end{document}